\documentclass[a4paper]{amsart}

\usepackage{amssymb,latexsym}
\usepackage{graphics}

\theoremstyle{plain}

\theoremstyle{remark}

\numberwithin{equation}{section}

\newcommand{\R}{\mathbb{R}}

\newcommand{\TV}{\text{\rm Tot.Var.}}

\newcommand{\la}{\bigl\langle}
\newcommand{\ra}{\bigr\rangle}

\begin{document}

\title[Singular limits to hyperbolic systems]{A note on singular limits to 
hyperbolic systems}

\author{Stefano Bianchini}
\address{S.I.S.S.A. - I.S.A.S., via Beirut 2-4, 34014 TRIESTE (ITALY)}
\email{bianchin@sissa.it}
\urladdr{http://www.sissa.it/\textasciitilde bianchin/}
\date{August 30, 2000}

\begin{abstract}
In this note we consider two different singular limits to hyperbolic system of 
conservation laws, namely the standard backward schemes for non linear
semigroups and the semidiscrete scheme.

Under the assumption that the rarefaction curve of the 
corresponding hyperbolic system are straight lines, we prove the 
stability of the solution and the convergence to the perturbed system to the 
unique solution of the limit system for initial data with small total
variation.  
\end{abstract}

\maketitle

\centerline{S.I.S.S.A. Ref. 85/2000/M}
\vskip 1cm

\section{Introduction}\label{S:intro}

Consider a hyperbolic system of conservation laws
\begin{equation}\label{E:hcl1}
\left\{ \begin{array}{ccc}
u_t + f(u)_x &=& 0 \\
u(0,x) &=& u_0(x)
\end{array} \right.
\end{equation}
where $u \in \R^n$ and $f$ is a smooth function from an open set
$\Omega \subseteq \R^n$ with values in $\R^n$. Let $K_0$ be a compact set 
contained in $\Omega$, and let $\delta_1$ sufficiently small such that the 
compact set 
\begin{equation}\label{E:compact1}
K_1 \doteq \Bigl\{ u \in \R^n: \text{dist}\bigl(u,K_0 \bigr) \leq \delta_1 \}
\end{equation}
is entirely contained in $\Omega$.

We assume that the Jacobian 
matrix $A = Df$ is uniformly strictly hyperbolic in $K_1$, i.e. 
\begin{equation}\label{E:strhyp1}
\min_{i<j} \Bigl\{ \lambda_j(u) - \lambda_i(u) \Bigr\} \geq c >0, \ \forall 
u,v \in K_1, 
\end{equation}
where we
denote by $\lambda_i$ the eigenvalues of $A$, $\lambda_i < \lambda_j$. Let 
$r_i$, $l^i$ be the its right, left eigenvectors. 

In this setting it is well known that if $u_0(-\infty) \in K_0$ and $\TV(u_0) 
$ is sufficiently small, there exists a unique "entropic" solution $u: 
[0,+\infty) \mapsto {u_0}+ L^1(\R, \R^n) \cap \text{BV}(\R,\R^n)$ in the sense 
of \cite{bre:unique}. Moreover these solutions can be constructed as limits of 
wave front tracking approximations and they depend Lipschitz continuously on 
the initial data.

For a special class of systems, called in \cite{brejen:godunov}
{\it Straight Line Systems}, i.e. systems such that 
\begin{equation}\label{E:straight1}
\bigl( Dr_i \bigr) r_i = 0,
\end{equation}
very recently it has been proved that solutions to~\eqref{E:hcl1} can
be constructed as $L^1$ limits of solutions to different singular approximation 
of the hyperbolic system:
\begin{itemize}
\item Vanishing viscosity approximation \cite{biabre:visco}. This is the limit as
$\epsilon \to 0$ of the solution $u^\epsilon(t)$ of the system
\[
u_t + f(u)_x - \epsilon u_{xx} = 0.
\]
\item Relaxation approximation \cite{bia:relax,breshe:chrom1}. While in case 1) the
perturbation is parabolic, in this case we consider a hyperbolic
perturbation, namely
\[
u_t + f(u)_x = \epsilon \bigl( \Lambda^2 u_{xx} - u_{tt} \bigr),
\]
where $\Lambda$ is strictly bigger than all the eigenvalues of
$Df(u)$. 
\item Godunov scheme \cite{brejen:godunov}. This is a discrete scheme obtained
from~\eqref{E:hcl1} by considering differential ratio instead of
derivatives:
\[
u(n+1,j+1) = u(n,j+1) + \frac{\Delta t}{\Delta x} \Bigl[ f\bigl(u(n,j)\bigr) -
f\bigl(u(n,j+1)\bigr) \Bigr],
\]
where for stability condition it is assumed that $0 < \lambda_1
<\dots<\lambda_n<\Delta x/ \Delta t$.
\end{itemize}
The main idea behind these approximations is to obtain uniform BV
estimates for $t \geq 0$, if the initial data $u_0$ are of
sufficiently small total variation. 

This task is achieved by 
decomposing the equations satisfied by $u_x$, or $u_x$ and $u_t$ in
\cite{bia:relax}, or $u(n,j) - u(n,j-1)$
in \cite{brejen:godunov}, as $n$ scalar perturbed
conservation laws, coupled by terms 
of higher order. These terms are then considered as the source 
of total variation. For the special case of straight line systems, a
decomposition of $u_x$ which makes the source terms integrable is the
projection along the eigenvectors $r_i$ of the Jacobian $Df(u)$:
\[
u_x = \sum_i v^i r_i.
\]
Once it is proved that the $L^1$ norm of the component $v^i$ is
bounded, by Helly's theorem there exists a subsequence
$u^{\epsilon_k}$ converging to a weak solution $\bar u(t)$ of~\eqref{E:hcl1} as $k
\to \infty$. 

To prove the uniqueness of the limit $\bar u(t)$, one consider the
equation for a perturbation $h$ of the singular approximations. 
We observe that  $h=u_x$ is a particular solution of such system.
A generalization of the arguments used to prove an a priori bound on the 
total variation of $u$ shows the boundedness of the
$L^1$ norm of the components $h^i$, where
\[
h = \sum_i h^i r_i(u).
\]
By a standard homotopy argument \cite{biabre:visco}, this yields the stability of
all solutions of the approximating system.  Since the Lipschitz continuous
dependence on the initial data is uniform w.r.t.~both $\epsilon$ and
$t$, in the limit we obtain a uniform Lipschitz semigroup.

Finally it is well known that a uniform Lipschitz semigroup of
solutions to~\eqref{E:hcl1} is uniquely defined if we know the jumps
conditions of the entropic shocks, see \cite{bre:libro}. In this case, because of the 
condition~\eqref{E:straight1} and because in the scalar case the 
solution $u^\epsilon$ converges to the entropic solution, an argument similar 
to the one in \cite{biabre:visco} implies that the jump conditions coincide 
with the scalar jumps along the eigenvectors $r_i$.

Thus, under the assumption~\eqref{E:straight1}, the limit semigroup is 
independent on the approximation and coincides with the solution constructed by 
wave front tracking using the classical Lax Riemann solver.

In this note we want to extend the previous approach to the following
cases:
\begin{enumerate}
\item semigroup approximation \cite{cra:semappr, ser:libro}. This is obtained as limit of 
the system 
\[
\frac{u(t,x) - u(t-\epsilon,x)}{\epsilon} + A \bigl( u(t,x) \bigr) u(t,x)_x = 
0.
\]
This is the standard backward scheme for non-linear semigroups.
\item Semi-discrete schemes \cite{ben:semi1}, i.e. infinite dimensional ODE defined by
\[
\frac{\partial}{\partial t}u(t,x) + \frac{1}{\epsilon}\Bigl( 
f\bigl(u(t,x)\bigr) - f\bigl(u(t,x-\epsilon)\bigr) \Bigr)= 0.
\]
\end{enumerate}
We will prove that as $\epsilon \to 0$ the limits of the respective
solutions converge to a unique solution to~\eqref{E:hcl1}, and that
this limits defines a Lipschitz continuous semigroup $\mathcal{S}$ on
the space of function with small TV. Moreover this semigroup is
perfectly defined by a Riemann solver which, as explained above, coincides 
with the classical one. 

The same can be proved for quasilinear systems as in 
\cite{biabre:visco,brejen:godunov}, but for simplicity we consider only systems 
on conservation forms. 

Without any loss of generality we assume that
\begin{equation}\label{E:genass1}
\min_i \bigl\{ \lambda_i(u) \bigr\} = \kappa >0, \qquad 
\max_i \bigl\{ \lambda_i(u) \bigr\} = K < 1,
\end{equation}
for all $u$ in the compact set $K_1$. The second condition is needed only in 
case 2).

\section{Approximation by semigroup theory}\label{S:apperox1}

We consider in this section the case 1) of Section~\ref{S:intro}, i.e. the 
following singular approximation to system of conservation laws:
\begin{equation}\label{E:singappr1}
\frac{u(t,x) - u(t-\epsilon,x)}{\epsilon} + A \bigl( u(t,x) \bigr) u(t,x)_x = 
0, 
\end{equation}
where we recall that $u \in \R^n$ and $A(u) = Df(u)$.
By the rescaling $t \to t/\epsilon$, $x \to x/ \epsilon$ and setting for 
simplicity $u_n(x) = u(n,x)$, we obtain the evolutionary equations
\begin{equation}\label{E:singappr2}
u_{n} - u_{n-1} + A \bigl( u_n \bigr) u_{n,x} = 0.
\end{equation}
It is easy to prove that if the BV norm of $u_{n-1}$ is sufficiently
small, then $u_n$ exists: in fact the solution can be represented as
\[
u_n(x) = \int_{-\infty}^x
\text{exp}\left\{\int_x^y A^{-1} \bigl( u_n(z)
\bigr) dz \right\} A^{-1} \bigl( u_n(y) \bigr) u_{n-1}(y) dy,
\]
and since the eigenvalues of $A$ are positive we have that 
\[
\bigl\| u_n \bigr\|_{\infty} \leq C \TV\bigl(u_{n-1} \bigr),
\]
$C$ being a uniform constant of $\|A^{-1}\|_{\infty}$ in the compact set $K_0$. 

\subsection{Projection on rarefaction curves}\label{SS:proj1}
 
We now start the procedure explained in Section \ref{S:intro}.
By projecting the derivative along the eigenvectors $r_i(u_n)$ of $A(u_n)$
\begin{equation}\label{E:proj1}
u_{n,x} = \sum_i v^i_n r_i(u_n) = \sum_i v^i_n r_{i,n},
\end{equation}
the equations for the components $v^i$ are
\[
\sum_i v^i_{n} r_{i,n} - \sum_i v^i_{n-1} r_{i,n-1}+ 
\sum_i \bigl( \lambda_{i,n}v^i_{n} r_{i,n} 
\bigr)_x = 0. 
\]
This can be rewritten as
\begin{equation}\label{E:compeq2}
\sum_i \Bigl(  v^i_{n} - v^i_{n-1} + \bigl( \lambda_{i,n} v^i_{n} \bigr)_x 
\Bigr) r_{i,n} = \sum_i v^i_{n-1} \bigl( r_{i,n-1} - r_{i,n} \bigr) - 
\sum_{i,j} \lambda_{i,n} v^i_{n} v^j_{n} \bigl( Dr_{i,n} \bigr) r_{j,n} 
\end{equation}
The left-hand side is in conservation form, and we consider the right-hand side 
as the source of total variation. If we assume as in the introduction that 
$(Dr_i ) r_i(u) = 0$, the function $r_{i}(u) - r_{i}(v)$ is zero when $u - v$ 
is parallel to $r_i (u) = r_i(v)$. Thus we have
\begin{equation}\label{E:expas1}
r_{i}(u) - r_{i}(v) = \sum_{j \not= i} \alpha_j(u,v) \la l^j(u), u - v \ra, 
\end{equation}
where $\alpha_j(u,u) = r_j(u)$. Using~\eqref{E:singappr2}, the 
expansion~\eqref{E:compeq2} thus becomes 
\begin{align}\label{E:compeq3}
v^i_{n} - v^i_{n-1} + \bigl( \lambda_{i,n} v^i_{n} \bigr)_x =&~ \sum_{j \not= 
k}  \Bigl( \lambda_{k,n} v^j_{n-1} v^k_n \la l^i_n, \alpha_j(u_n,u_{n-1}) \ra - 
v^i_{n} v^j_{n} \la l^i_n, \bigl( Dr_{i,n} \bigr) r_{j,n} \ra \Bigr) \\
=&~ \sum_{j \not= k} H^i_{jk}(n) v^j_{n-1} v^k_n + \sum_{j \not= k} K^i_{jk}(n) 
v^j_{n} v^k_n. \notag
\end{align}
To estimate the source terms in~\eqref{E:compeq3}, we first consider the case 
of two linear equations.

\subsection{Analysis of the linear case}\label{S:linear}

\begin{figure}
\centerline{\resizebox{16cm}{6cm}{{\includegraphics{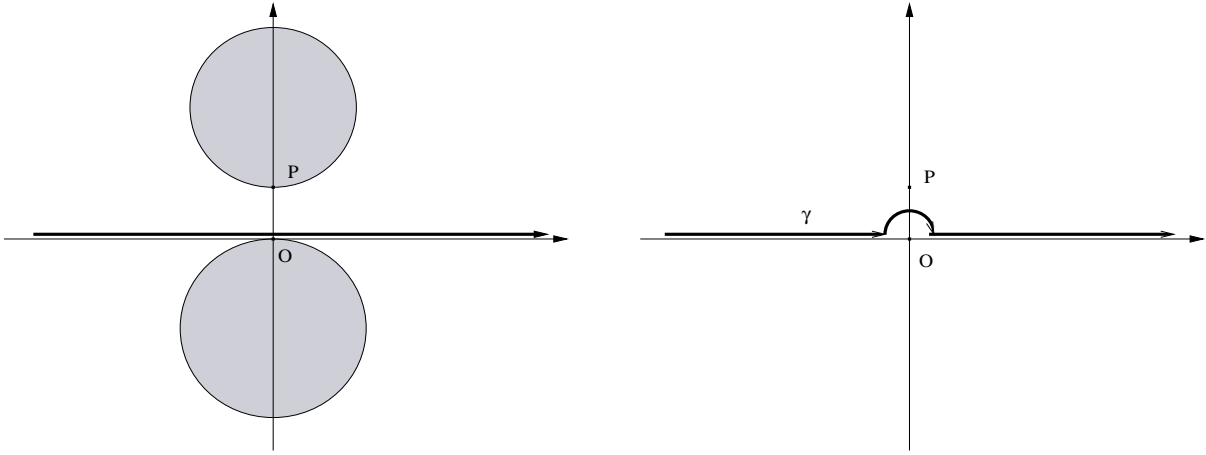}}}}
\caption{Integration path on the complex plane, where $P = i (\lambda
- \mu)/\lambda \mu$.} 
\label{Fi:integ1}
\end{figure}

Consider a single linear equation
\begin{equation}\label{E:linear1}
v_n - v_{n-1} + \lambda v_{n,x} = 0, \qquad \lambda > 0.
\end{equation}
We can find the fundamental solution to the previous equation by means of 
Fourier transform: we have
\[
v_n(x) = \int_{\R} c(n,\xi) e^{-i\xi x} d\xi,
\]
and substituting
\[
c(n,\xi) - c(n-1,\xi) - i \lambda \xi c(n,\xi) = 0 \quad \Longrightarrow 
\quad c(n,\xi) = \frac{c_0(\xi)}{\bigl( 1 - i \lambda \xi \bigr)^n}.
\]
In particular the fundamental solution has $c_0 (\xi) \equiv 1 / 2 \pi$, so 
that 
\begin{equation}\label{E:fundsol1}
v_n(x) = 
\frac{1}{\lambda} \left( \frac{x}{\lambda} \right)^{n-1} \frac{e^{-x / 
\lambda}}{\bigl(n-1\bigr)!} \chi_{[0,+\infty)}(x).
\end{equation}
Consider two equations of the form~\eqref{E:linear1},
\begin{align}\label{E:linear2}
v_n - v_{n-1} + \lambda v_{n,x} = 0 \\
z_n - z_{n-1} + \mu z_{n,x} = 0 \notag
\end{align}
with initial data $v_0(x) = \delta(x)$ and $z_0(x) = \delta(x - x_0)$,
and assume without any loss of generality that $\lambda > \mu >0$.
We can compute the intersection integrals: denoting with $d(n,\xi)$ the Fourier 
coefficients of $z_n(x)$ we have 
\begin{align}\label{E:transcomp1}
\sum_{n=0}^{N} \int_{\R} v_n(x) z_n(x) dx =&~ \sum_{n=0}^{N} 
2 \pi \int_{\R} c(n,\xi) d(n,-\xi) e^{-i \xi x_0} d\xi \\
=&~ \frac{1}{2 \pi} \int_{\R} \sum_{n=0}^{N} \frac{1}{\bigl( 1 - i 
\lambda \xi \bigr)^n  \bigl( 1 + i \mu \xi \bigr)^n} e^{-i \xi x_0} d\xi. \notag 
\end{align}
If $\xi$ is considered as a complex variable,
we can let $N \to +\infty$ only in the region where 
\[
Z \doteq \Bigl\{ \xi \in \mathbb{Z}: \bigl|( 1 - i 
\lambda \xi )( 1 + i \mu \xi )| < 1 \Bigr\}
\]
i.e. outside the regions depicted
in Figure \ref{Fi:integ1}. Deforming the path to avoid the region $Z$, we can 
pass to the limit: 
\begin{align}
\sum_{n=0}^{+\infty} \int_{\R} v_n(x) z_n(x) dx =&~ \frac{1}{2 \pi}
\int_{\gamma} \frac{e^{-i \xi x_0}}{1-\frac{1}{( 1 - i 
\lambda \xi ) ( 1 + i \mu \xi )}} d \xi \notag \\
=&~ \frac{1}{2 \pi} \int_{\gamma} \frac{( 1 - i 
\lambda \xi ) ( 1 + i \mu \xi )}{( 1 - i 
\lambda \xi ) ( 1 + i \mu \xi ) - 1} e^{-i \xi x_0} d\xi. \notag
\end{align}
By means of complex analysis we have finally that 
\begin{equation}\label{E:transcomp2}
P(x_0) \doteq \sum_{n=0}^{+\infty} \int_{\R} v_n(x) z_n(x) dx = \begin{cases}
1/(\lambda - \mu) \cdot \text{exp}\Bigl( (\lambda - \mu)/(\lambda \mu) x_0 
\Bigr) & x_0 < 0 \\
1/(\lambda - \mu) & x_0 \geq 0 
\end{cases}
\end{equation}

\noindent In fact, depending on the sign of $x_0$, the integration along the
line $\gamma$ is equivalent to the integration around the pole $0$ or the 
pole $P=i (\lambda - \mu)/\lambda \mu$.

\subsection{BV estimates}\label{S:bvest}

Now to prove that~\eqref{E:singappr2} has a solution with uniformly bounded 
total variation. Define the functional 
\begin{equation}\label{E:interpot1}
Q(n) = Q(u_n,u_{n-1}) \doteq \sum_{i < j} \int_{\R} P_0(x - y) \Bigl\{ 
\bigl| v^i_n(x) v^j_n(y) \bigr| + \bigl| v^i_{n-1}(x) v^j_n(y) \bigr| + 
\bigl| v^i_{n}(x) v^j_{n-1}(y) \bigr| \Bigr\} dx dy , 
\end{equation}
where $P$ is computed substituting to $\lambda-\mu$ the constant of separation 
of speeds $c$, and taking the minimal value of the exponent $(\lambda - 
\mu)/\lambda \mu$: 
\[
P_0(x) \doteq 
\begin{cases}
1/c \cdot \text{exp}\Bigl( c/\bigl(K (K-c)\bigr) x_0 
\Bigr) & x_0 < 0 \\
1/c & x_0 \geq 0 
\end{cases}
\]
We recall that $c$ and $Kappa$ are defined in the introduction.

Using the same analysis of \cite{bia:relax}, we see immediately 
that 
\begin{align}\label{E:derivpot1}
Q(n) - Q(n-1) =&~ \sum_{i < j} \int_{\R} P_0(x - y) \Bigl\{ \bigl| v^i_n(x) 
v^j_n(y) \bigr| - \bigl| v^i_{n-1}v^j_{n-1} \bigr| \Bigr\} \\
& \qquad + \sum_{i<j} \int_{\R} P_0(x-y) \Bigl\{ \bigl| v^i_{n-1}(x) v^j_n(y) 
\bigr| - \bigl| v^i_{n-2}v^j_{n-1} \bigr| \Bigr\} \notag \\
& \qquad + \sum_{i<j} \int_{\R} P_0(x-y) \Bigl\{ \bigl| v^i_{n}(x) v^j_{n-1}(y) 
\bigr| - \bigl| v^i_{n-1}v^j_{n-2} \bigr| \Bigr\} dx dy \notag \\
\leq&~ - \Bigl(1 - C \max_{m=1, \dots, n} \TV \bigl(u_m \bigr) \Bigr)\left[ 
\sum_{j \not= k} \bigl| v^j_{n-1} v^k_n \bigr| + \sum_{j \not= k} \bigl| v^j_n 
v^k_n \bigr| \right], \notag
\end{align}
where $C$ is c constant depending only on $H^i_{jk}$, $K^i_{jk}$, $\kappa$, $K$ 
and $c$. Thus if $\delta_0$ is sufficiently small, using~\eqref{E:derivpot1} we 
have 
\[
\TV(u_1) + C_0 Q(u_1,u_0) \leq \delta_1 \quad \text{and} \quad
\frac{d}{dt} \Bigl\{ \TV(u) + C_0 Q(u) \Bigr\} \leq 0,
\]
where the constant $C_0$ is big enough, independent on $\delta_0$. 
This proves that the solution $u_n$ has uniformly bounded total
variation for all $n \in \mathbb{N}$.

\subsection{Stability estimates}\label{S:stabil}

We now consider the stability estimates of~\eqref{E:singappr2}. The 
equations for a perturbation $u + \delta h$ as $\delta \to 0$ are
\begin{equation}\label{E:pertueq1}
h_n - h_{n-1} + \bigl( A(u_n) h_n \bigr)_x = \bigl( DA\bigl(u_n\bigr) u_{n,x} 
\bigr) h - \bigl( DA\bigl(u_n\bigr) h_n \bigr) u_{n,x}.
\end{equation}
Using the same projection of~\eqref{E:proj1}, i.e. 
\[
h_n = \sum_i h^i_n r_{i,n},
\]
we have that the equations for the components $h^i_n$ are
\begin{align}\label{E:perteq2}
h^i_{n} - h^i_{n-1} + \bigl( \lambda_{i,n} h^i_{n} \bigr)_x =&~ \sum_{j \not= 
k}  \Bigl( \lambda_{k,n} h^j_{n-1} v^k_n \la l^i_n, \alpha_j(u_n,u_{n-1}) \ra - 
h^i_{n} v^j_{n} \la l^i_n, \bigl( Dr_{i,n} \bigr)r_{j,n}\ra \Bigr) \\
& \qquad \qquad \qquad + \sum_{j \not= k} h^i_n v^j_n 
\la l^i_n,  \bigl( A(u_n) r_{j,n} \bigr) r_{i,n} - \bigl( A(u_n) r_{i,n} \bigr)
r_{j,n} \ra  \notag \\ 
=&~  \sum_{j \not= k} H(n) h^j_{n-1} v^k_n + \sum_{j \not= k} K'(n) 
h^j_{n} v^k_n. \notag 
\end{align}
Using the same analysis of \cite{bia:relax}, it is easy to prove that
a functional as in Section~\ref{S:bvest} gives the stability of the
solution. 

\section{Approximation by semi-discrete scheme}\label{S:apperox2}

We now consider the case 2) of Section~\ref{S:intro}, i.e. the following 
singular approximation to system of conservation laws:
\begin{equation}\label{E:semidiscr1}
\frac{\partial}{\partial t}u(t,x) + \frac{1}{\epsilon}\Bigl( 
f\bigl(u(t,x)\bigr) - f\bigl(u(t,x-\epsilon)\bigr) \Bigr)= 0, 
\end{equation}
where $u \in \R^n$.
By the rescaling $t \to t/\epsilon$, $x \to x/ \epsilon$, we obtain the 
evolutionary equations
\begin{equation}\label{E:semidiscr2}
\dot u_{n}(t) +  f\bigl(u_n(t)\bigr) - f\bigl(u_{n-1}(t)\bigr)= 0.
\end{equation}
The equation for the "derivative" $v_n \doteq u_n - u_{n-1}$ are
\begin{equation}\label{E:deriveq1}
\dot v_{n}(t) +  f\bigl(u_n(t)\bigr) - 2 f\bigl(u_{n-1}(t)\bigr) + 
f\bigl(u_{n-2}(t)\bigr)= 0. 
\end{equation}

\subsection{Projection on rarefaction curves}\label{SS:proj2}

The vector $v_n$ is now decomposed along the eigenvectors 
$r_{i,n}$ of the Riemann problem $u_{n-1}$, $u_n$: we have
\[
\dot u_n(t) + \sum_i \lambda_{i,n} v^i_n r_{i,n} = 0,
\]
\begin{align}\label{E:deriveq2}
\sum_i \Bigl( \dot v^i_n + \lambda_{i,n} v^i_n - \lambda_{i,n-1} 
v^i_{n-1} \Bigr) r_{i,n} =&~ - \sum_{i,j} v^i_n v^j_n \bigl(D r_{i,n}\bigr) 
r_{j,n} -  \sum_ {i,j} v^i_n v^j_{n-1} \bigl( Dr_{i,n} \bigr) r_{j,n-1} \\
& \qquad \qquad + \sum_i \lambda_{i,n-1} v^i_{n-1} \bigl( r_{i,n-1} - r_{i,n} 
\bigr) \notag \\
=&~ - \sum_{i,j} v^i_n v^j_n  \bigl(D r_{i,n}\bigr) r_{j,n}-  
\sum_ {i,j} v^i_n v^j_{n-1} \bigl( Dr_{i,n} \bigr) r_{j,n} \notag \\
& \qquad \qquad + \sum_ {i,j} v^i_n v^j_{n-1} \bigl( Dr_{i,n} \bigr)\bigl( 
r_{j,n} - r_{j,n-1} \bigr) \notag \\
& \qquad \qquad  + \sum_i \lambda_{i,n-1} v^i_{n-1} \bigl( 
r_{i,n-1} - r_{i,n}  \bigr),  \notag 
\end{align}
where $\lambda_{i,n}$ and $r_{i,n}$ are the eigenvalues and right eigenvectors 
of the average matrix
\[
A\bigl(u_n,u_{n-1}\bigr) \doteq \int_0^1 Df\Bigl( u_{n-1} + \bigl( u_n - 
u_{n-1} \bigr) s \Bigr) ds.
\]

If we assume the condition~\eqref{E:straight1}, 
the functions $(Dr_{i,n})r_{j,n}$ and $r_{i,n} - r_{i,n-1}$ are zero when 
$u_n - u_{n-1}$ and $u_{n-1} - u_{n-2}$ are parallel to $r_{i,n} = r_{i,n-1}$. 
Thus we have 
\begin{align}\label{E:expans1}
r_{j,n}\bullet r_{i,n} =&~ \sum_{j \not= i}\alpha_{j,n} v^j_n, \\
r_{i,n} - r_{i,n-1} =&~ \sum_{j \not= i} \beta_{j,n} v^j_n + \sum_{j \not= i} 
\gamma_{j,n-1} v^j_{n-1}, \notag
\end{align}
as in Section~\ref{SS:proj1}.
Using~\eqref{E:expans1}, the 
expansion~\eqref{E:deriveq2} thus becomes 
\begin{equation}\label{E:deriveq3}
\dot v^i_{n} + \lambda_{i,n} v^i_n - \lambda_{i,n-1} 
v^i_{n-1} = \sum_{j \not= k} H_n(t) v^j_n v^k_n + \sum_{j \not= k} G_n(t)
v^j_n  v^k_{n-1}. \end{equation}
To estimate the source terms in~\eqref{E:deriveq3}, we consider the case of two 
linear equations.

\subsection{Analysis of the linear case}\label{S:linear2}

\begin{figure}
\centerline{\resizebox{16cm}{6cm}{{\includegraphics{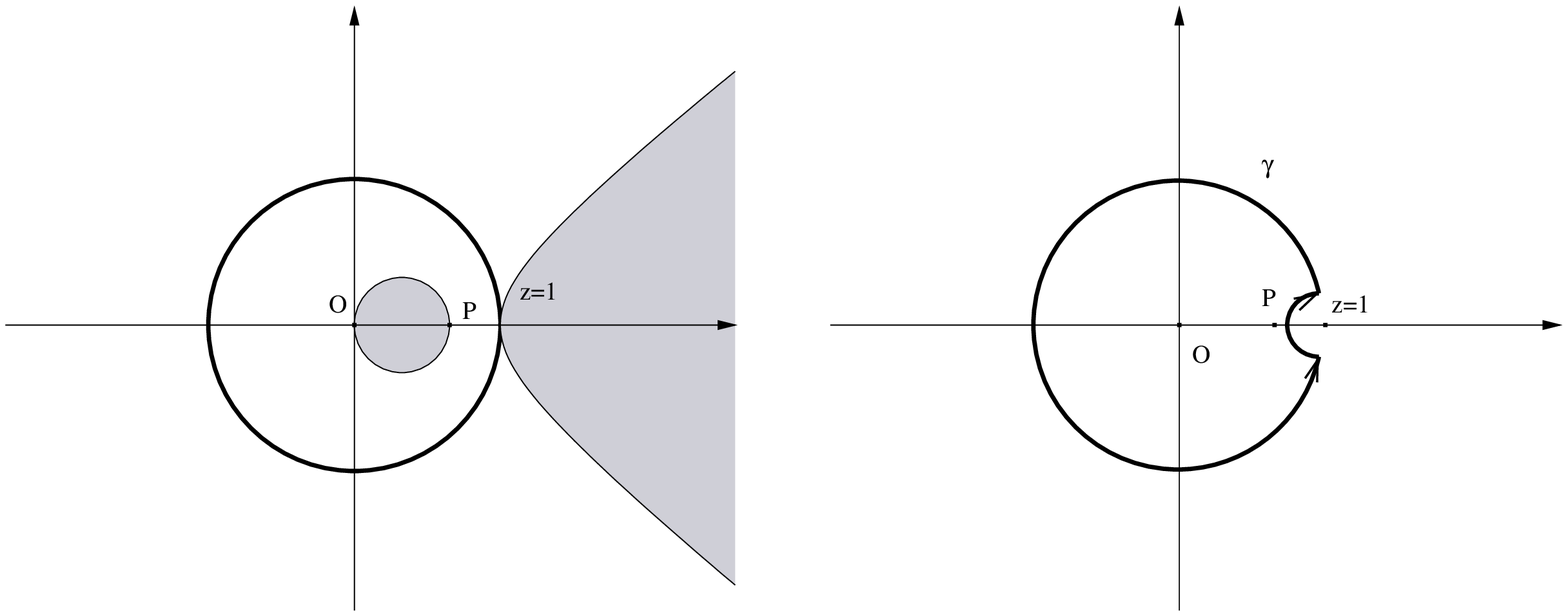}}}}
\caption{Integration path on the complex plane, where $P = \mu/\lambda$.} 
\label{Fi:integ2}
\end{figure}

Consider a single linear equation
\begin{equation}\label{E:linear3}
\dot v^i_{n} + \lambda v^i_n - \lambda 
v^i_{n-1} = 0, \qquad \lambda > 0.
\end{equation}
We can find the fundamental solution to the previous equation by means of 
Fourier transform: defining the periodic function
\[
c(t,x) \doteq \sum_n v_n(t) e^{i n x},
\]
we have that the equation satisfied by $c$ is
\begin{align}
c_t =&~ \sum_n \dot v_n e^{i n x} = \lambda \sum_n \bigl( v_{n-1} - v^n 
\bigr) e^{i n x} \notag \\
=&~ \lambda \bigl( e^{ix} - 1 \bigr) c, \notag
\end{align}
whose general solution is 
\[
c(t,x) = c(0,x) \text{exp} \Bigl( \lambda \bigl( e^{ix} - 1 \bigr) t \Bigr).
\]
In particular the fundamental solution starting at $n_0$ has $c(0,x) = 
\text{exp}(i n_0 x)$, so that if $n_0 = 0$
\begin{equation}\label{E:fundsol2}
v_n(t) = \begin{cases}
0 & x < 0 \\
\bigl( \lambda t \bigr)^n / n! \cdot \text{exp}\bigl(- \lambda t \bigr) & x 
\geq 0 
\end{cases}
\end{equation}
If now we consider two equations of the form~\eqref{E:linear3},
\begin{align}\label{E:linear4}
\dot v_n + \lambda \bigl( v_{n} - v_{n-1} \bigr) = 0, \\
\dot z_n + \mu \bigl( z_{n} - z_{n-1} \bigr) = 0. \notag
\end{align}
we can compute the intersection integrals: denoting with $d(t,x)$ the Fourier 
transform of $z_n(t)$ and assuming that $\lambda > \mu >0$, we have 
\begin{align}\label{E:transcomp3}
\int_0^{+\infty} \sum_{n=-\infty}^{+\infty} v_n(t) z_n(t) dt =&~ 
\int_{0}^{+\infty} \frac{1}{2 \pi} \int_{0}^{2 \pi} c(t,x) d(t,-x) e^{-i n_0 x} 
dx dt\\ 
=&~ \frac{1}{2 \pi} \int_{0}^{2 \pi} \int_{0}^{+\infty} \text{exp} \Bigl( 
\lambda \bigl( e^{ix} - 1 \bigr) t + \mu \bigl( e^{-ix} - 1 
\bigr) t \Bigr) e^{-i n_0 x} dx dt \notag \\ 
=&~ \frac{1}{2 \pi} \int_{0}^{2 \pi}
\frac{e^{-i n_0 x}}{\lambda \bigl( e^{ix} - 1 \bigr) + \mu \bigl( e^{-ix} - 1 
\bigr) } dx \notag \\ 
=&~ \frac{1}{2 \pi i} \oint_{\gamma}
\frac{z^{- n_0}}{\bigl( z - 1 \bigr) \bigl( \lambda z - \mu \bigr) } dz,
\notag 
\end{align}
where $\gamma$ is the path represented in Figure~\ref{Fi:integ2}. 

By means of complex analysis we have that 
\begin{equation}\label{E:transcomp4}
P(n_0) \doteq \int_{0}^{+\infty} \sum_{n=-\infty}^{+\infty}  v_n(t) z_n(t) dx = 
\begin{cases} 
1/(\lambda - \mu) \cdot \bigl( \lambda / \mu \bigr)^{n_0} & n_0 < 0 \\
1/(\lambda - \mu) & n_0 \geq 0 
\end{cases}
\end{equation}

\subsection{BV estimates}\label{S:bvest2}

Now to prove that~\eqref{E:semidiscr1} has a solution with uniformly bounded 
total variation. By defining the functional 
\begin{equation}\label{E:interpot2}
Q\bigl( u(t) \bigr) \doteq \sum_{i < j} \sum_{n,m=-\infty}^{+\infty} P(n - m) 
\Bigl\{ v^i_n(t) v^j_m(t) + v^i_{n-1}(t) v^j_m(t) + v^i_{n}(t) v^j_{m-1}(t) 
\Bigr\} dx dy , 
\end{equation}
where $P$ is computed using the constant of separation of speeds $c$ instead of 
$\lambda - \mu$ and $1+c/K$ instead of $\lambda/\mu$, since the left hand side 
of~\eqref{E:deriveq3} is in conservation form, we conclude immediately that 
\[
\TV\bigl(u(0)\bigr) + C_0 Q\bigl(u(o)\bigr) \leq \delta_1, \qquad \frac{d}{dt} 
\Bigl\{ \TV(u) + C_0 Q(u) \Bigr\} \leq 0. 
\]
This concludes the proof of bounded total variation.

\subsection{Stability estimates}\label{S:stabil2}

Finally we consider the stability estimates of~\eqref{E:semidiscr1}. The 
equations for a perturbation $u + \delta h$ as $\delta \to 0$ are
\begin{equation}\label{E:pertueq3}
\dot h_n(t) + Df \bigl( u_n \bigr) h_n - Df \bigl( u_{n-1} \bigr) = 0.
\end{equation}
Considering the projection
\[
h_n(t) = \sum_i h^i_n(t) r_{i}(u_n),
\]
we have that the equations for the components $h^i_n$ are
\[
\dot h^i_{n} + \lambda_{i}(u_n) h^i_{n} - \lambda_{i}(u_{n-1}) h^i_{n-1} =
\sum_{j \not= k} H'(n) h^j_{n-1} v^k_n + \sum_{j \not= k} G'(n) 
h^j_{n} v^k_n. 
\]
At this point it is clear that a functional as in Section~\ref{S:bvest} proves 
the stability of the solution. This concludes the proof.

\vskip 0.2cm

\noindent{\bf Acknowledgment.} This research was
partially supported by the European TMR Network on Hyperbolic
Conservation Laws ERBFMRXCT960033.


\end{document}